\documentclass{article}
\usepackage{amssymb}
\usepackage{amsfonts}
\usepackage{amsmath}

\newtheorem{theorem}{Theorem}

\newtheorem{definition}[theorem]{Definition}

\begin{document}

\title{Four-dimensional Conformally Flat\\
Berwald and Landsberg Spaces }
\author{Gauree Shanker \\
Centre for Mathematics and Statistics\\
Central University of Punjab, Bathinda-151001\\
Email:gshankar@cup.ac.in}
\date{}
\maketitle

\begin{abstract}
{\small {{\footnotesize The problem of conformal transformation and
conformal flatness of Finsler spaces has been studied by so many researchers 
$\left[ 6,16,17,20,21\right] .$ Recently, Prasad et. al $\left[ 19\right] $
have studied three dimensional conformally flat Landsberg and Berwald spaces
and have given some important results. The purpose of the present paper is
to extend the idea of conformal change to four dimensional Finsler spaces
and find the suitable conditions under which a four dimensional conformally
flat Landsberg space becomes a Berwald space. } \vspace{0.2in}\newline
}}
\end{abstract}

\begin{quotation}

\end{quotation}

\section{\textbf{Introduction }{\protect\small {\protect\vspace{0.2in}}}}

{\small {Let us consider a four-dimensional Finsler space $F^4=(M,L)$
equipped with the fundamental function $L\left( x,y\right) .$ The
fundamental metric tensor $g_{ij}$ and Cartan $C-$tensor $C_{ijk}$ of $F^{4}$
are defined by ([2], [14]) }}

\begin{description}
\item[$\left( 1.1\right) $] {\small $g_{ij}=\frac{1}{2}\dot{\partial }_{i}%
\dot{\partial}_{j}L^{2},\qquad C_{ijk}=\frac{1}{2}\dot{\partial }_{k}g_{ij}=%
\frac{1}{4}\dot{\partial }_{i}\dot{\partial }_{j}\dot{\partial }_{k}L^{2}.$ }
\end{description}

{\small Throughout the paper, the symbols $\dot{\partial}_{i}=\frac{\partial 
}{\partial y^{i}}$ and $\partial _{i}=\frac{\partial }{\partial x^{i}}$ have
been used. Let $\left( x,y\right) =\left( x^{i},y^{i}\right) $ be a local
coordinate system of the total space of the tangent bundle }${\small TM}$%
{\small \ of the four- dimensional differentiable manifold }${\small M}$%
{\small . If in a Finsler space there exists a local coordinate system $%
\left( x^{i},y^{i}\right) $ in which the fundamental tensor $g_{ij}$ can be
written as a function of the variables $y^{i}$ alone, the space is called
locally Minkowski space and such a coordinate system is called a rectilinear
coordinate system. If a Finsler space $(M,L)$ is conformal to a locally
Minkowski space ($M,\overline{L}),$ then $(M,L)$ is called a conformally
flat Finsler space, more precisely, $\sigma -$ conformally flat Finsler
space.}

\section{\textbf{Scalar components in Miron frame}}

{\small Let $F^{4}=(M, L)$ be a {\normalsize four-dimensional Finsler space
with fundamental function }$L\left( x,y\right) $. The metric tensor $g_{ij}$
and Cartan $C-$tensor $C_{ijk}$ of $F^{4}$ are defined by (1.1). The frame $%
\left\{ e_{(\alpha )}^{i}\right\} $ $(\alpha =1,...,4)$ is called the Miron
frame of $F^{4}$, where $e_{\left( 1\right) }^{i}=l^{i}=\frac{y^{i}}{L}$ is
called the normalized supporting element, $e_{\left( 2\right) }^{i}=m^{i}=%
\frac{C^{i}}{C}$ is called the normalized torsion vector, $e_{\left(
3\right) }^{i}=n^{i},e_{\left( 4\right) }^{i}=p^{i}$ are constructed by $%
g_{ij}e_{\left( \alpha \right) }^{i}e_{\left( \beta \right) }^{j}=\delta
_{\alpha \beta }.$ Here, $C$ is the length of torsion vector $%
C_{i}=C_{ijk}g^{jk}.$ The Greek letters $\alpha ,\beta ,\gamma ,\delta $
vary from $1$ to }${\small 4}${\small \ throughout the paper. Summation
convention is applied for both the Greek and Latin indices. In the Miron
frame an arbitrary tensor can be expressed by scalar components along the
unit vectors $l^{i},m^{i},n^{i},p^{i}.$ For instance; let $T=T_{j}^{i}$ be a
tensor field of $\left( 1,1\right) $ type, then the scalar components $%
T_{\alpha \beta }$ of $T$ are defined by $T_{\alpha \beta
}=T_{j}^{i}e_{\left( \alpha \right) i}e_{\left( \beta \right) }^{j}$and the
components $T_{j}^{i}$ of the tensor $T$ are expressed as $%
T_{j}^{i}=T_{\alpha \beta }e_{\left( \alpha \right) }^{i}e_{\left( \beta
\right) j}.$\newline
It is well known (See [18], [20]) that in a four-dimensional Finsler space
there are eight main scalars $H,\ I,\ J,\ K,\ H^{\prime },\ I^{\prime },\
J^{\prime },\ K^{\prime }$ in which the sum of $H,\ I$ and $K$ is $LC,$
which is called unified main scalar. In [20], it has been proved also that
in a four-dimensional Finsler space there exist three h-connection vectors $%
h_{i},j_{i},k_{i},\ $ whose scalar components with respect to the frame $%
\left\{ e_{\left( \alpha \right) }^{i}\right\} $ are $h_{\gamma },j_{\gamma
},k_{\gamma }$ i.e., $h_{i}=h_{\gamma }e_{\left( \gamma \right) i},\quad
j_{i}=j_{\gamma }e_{\left( \gamma \right) i}, k_{i}=k_{\gamma }e_{\left(
\gamma \right) i}$. Further, in a four-dimensional Finsler space there exist
three v-connection vectors $u_{i},v_{i},w_{i}$ whose scalar components with
respect to the frame $\left\{ e_{\left( \alpha \right) }^{i}\right\} $ are $%
u_{\gamma },v_{\gamma },w_{\gamma }\ i.e.,$ } {\small $u_{i}=u_{\gamma
}e_{\left( \gamma \right) i},\qquad v_{i}=v_{\gamma }e_{\left( \gamma
\right) i},\qquad w_{i}=w_{\gamma }e_{\left( \gamma \right) i}.$ }

{\small In the next sections, we use the following important results,
obtained in $\left[ 20\right]$: }

\begin{description}
\item[$\left( 2.1\right) \left( a\right) $] {\small $l_{i\shortmid j}=0,$ $%
m_{i\shortmid j}=n_{i}h_{j}+p_{i}j_{j},$ $n_{i\shortmid
j}=-m_{i}h_{j}+p_{i}k_{j},$ $p_{i\shortmid j}=-m_{i}j_{j}-n_{i}k_{j},$ }
\end{description}

{\small and }

\begin{description}
\item[$\left( 2.1\right) \left( b\right) $] {\small $Ll_{i\mid j}=$ $%
m_{i}m_{j}+n_{i}n_{j}+p_{i}p_{j}=g_{ij}-l_{i}l_{j}=h_{ij},$ }

\item {\small $\qquad Lm_{i\mid j}=-l_{i}m_{j}+n_{i}u_{j}+p_{i}v_{j},\qquad
Ln_{i\mid j}=-l_{i}n_{j}-m_{i}u_{j}+p_{i}w_{j},$ }

\item {\small $\qquad Lp_{i\mid j}=-l_{i}p_{j}-m_{i}v_{j}-n_{i}w_{j}.$ }
\end{description}

\section{\textbf{Conformal change of Cartan's connection}}

{\ Let $\left( M,L\right) $ and $\left( M,\bar{L}\right) $ be two Finsler
manifolds. Then the two manifolds are conformal if there exists a positive
differentiable function $\sigma :M\longrightarrow 
\mathbb{R}
$ such that $\bar{L}\left( x,y\right) =e^{\sigma \left( x\right) }L\left(
x,y\right) $. In this case, the transformation $L\longrightarrow \bar{L}$ is
said to be conformal transformation and the Finsler manifolds $\left(
M,L\right) $ and $\left( M,\bar{L}\right) $ are said to be conformal or
conformally related. If $\sigma \left( x\right) =$ constant, the
transformation is called homothetic, otherwise, it is called non-homothetic
transformation.\newline
We consider a conformal change $L\left( x,y\right) \longrightarrow \overline{%
L}\left( x,y\right) =e^{\sigma \left( x\right) }L\left( x,y\right) $ of a
four-dimensional Finsler space $F^{4}=(M^{4},L(x,y))$ with the fundamental
function $L(x,y)$, where $\sigma (x)$ is a scalar function of position $x^{i}
$ alone, called the conformal factor. We shall denote the Finsler space with
changed fundamental function $\overline{L}\left( x,y\right) $ by $\overline{F%
}^{4}=(M^{4},\ \overline{L}(x,y))$ and quantities of $\overline{F}^{4}$ by
upper line. The following change of important quantities are well known $%
([6])$. }

\begin{description}
\item {\small \ }

\item[$\left( 3.1\right) $] {\small $\overline{l}_{i}=e^{\sigma }l_{i},$ $%
\overline{m }_{i}=e^{\sigma }m_{i},$ $\overline{n}_{i}=e^{\sigma }n_{i},$ $%
\overline{p} _{i}=e^{\sigma }p_{i},$ $\overline{g}_{ij}=e^{2\sigma }g_{ij},$ 
}

\item[$\left( 3.2\right) $] {\small $\overline{l}^{i}=e^{-\sigma }l^{i},$ $%
\overline{ m}^{i}=e^{-\sigma }m^{i},$ $\overline{n}^{i}=e^{-\sigma }n^{i},$ $%
\overline{p }^{i}=e^{-\sigma }p^{i},$ $\overline{g}^{ij}=e^{-2\sigma
}g^{ij}, $ }

\item[$\left( 3.3\right) $] {\small $\overline{C}_{ijk}=e^{2\sigma }C_{ijk},$
$\overline{C}_{jk}^{i}=C_{jk}^{i},$ $\overline{H}=H,$ $\overline{I}=I,$ $%
\overline{J}=J,$ $\overline{K}=K,$ $\overline{H^{\prime }}=H^{\prime },$ $%
\overline{I^{\prime }}=I^{\prime },$ $\overline{J^{\prime }}=J^{\prime },$ $%
\overline{K^{\prime }}=K^{\prime }.$\qquad \qquad }
\end{description}

In [20], it has been obtained that:

\begin{description}
\item[$\left( 3.4\right) $] {\small $\overline{G}_{j}^{i}=G_{j}^{i}+Ll^{i}%
\left( \sigma _{1}l_{j}+\sigma _{2}m_{j}+\sigma _{3}n_{j}+\sigma
_{4}p_{j}\right) -Lm^{i}\left( \sigma _{2}l_{j}-\sigma _{5}m_{j}+\sigma
_{6}n_{j}+\sigma _{7}p_{j}\right) $ }

\item {\small $\qquad \qquad -Ln^{i}\left( \sigma _{3}l_{j}+\sigma
_{6}m_{j}-\sigma _{8}n_{j}-\sigma _{9}p_{j}\right) -Lp^{i}\left( \sigma
_{4}l_{j}+\sigma _{7}m_{j}-\sigma _{9}n_{j}-\sigma _{10}p_{j}\right) .$ }
\end{description}

{\small For the conformal change of the adapted components $h_{\alpha
},j_{\alpha },k_{\alpha }$ of three h-connection vectors $h_{i},j_{i},k_{i}$%
, from }${\small (3.1)}${\small \ and }${\small (2.1)(a)}${\small , we have
[20]}

\begin{description}
\item {\small \ }

\item[$\left( 3.5\right) \left( a\right) $] {\small $\overline{h}%
_{j}=h_{j}+\left\{ \sigma _{2}u_{2}+\sigma _{3}u_{3}+\sigma
_{4}u_{4}\right\} l_{j}$ }

\item {\small $\qquad \qquad +\left\{ -\sigma _{5}u_{2}+\sigma
_{6}u_{3}+\sigma _{7}u_{4}+\sigma _{5}\left( J+J^{\prime }\right) +\sigma
_{6}I+\sigma _{7}K^{\prime }-\sigma _{12}\right\} m_{j}$ }

\item {\small $\qquad \qquad +\left\{ \sigma _{6}u_{2}-\sigma
_{8}u_{3}-\sigma _{9}u_{4}-\sigma _{6}\left( J+J^{\prime }\right) -\sigma
_{8}I-\sigma _{9}K^{\prime }+\sigma _{14}\right\} n_{j}$ }

\item {\small $\qquad \qquad +\left\{ \sigma _{7}u_{2}-\sigma
_{9}u_{3}-\sigma _{10}u_{4}-\sigma _{7}\left( J+J^{\prime }\right) -\sigma
_{9}I-\sigma _{10}K^{\prime }+\sigma _{15}\right\} p_{j},$ }

\begin{description}
\item {\small \ }

\item[$\left( b\right) $] {\small $\overline{j}_{j}=j_{j}+\left\{ \sigma
_{2}v_{2}+\sigma _{3}v_{3}+\sigma _{4}v_{4}\right\} l_{j}$ }
\end{description}

\item {\small $\qquad \qquad +\left\{ -\sigma _{5}v_{2}+\sigma
_{6}v_{3}+\sigma _{7}v_{4}+\sigma _{5}\left( H^{\prime }+I^{\prime }\right)
+\sigma _{6}K^{\prime }+\sigma _{7}K-\sigma _{13}\right\} m_{j}$ }

\item {\small $\qquad \qquad +\left\{ \sigma _{6}v_{2}-\sigma
_{8}v_{3}-\sigma _{9}v_{4}-\sigma _{6}\left( H^{\prime }+I^{\prime }\right)
-\sigma _{8}K^{\prime }-\sigma _{9}K+\sigma _{15}\right\} n_{j}$ }

\item {\small $\qquad \qquad +\left\{ \sigma _{7}v_{2}-\sigma
_{9}v_{3}-\sigma _{10}v_{4}-\sigma _{7}\left( H^{\prime }+I^{\prime }\right)
-\sigma _{9}K^{\prime }-\sigma _{10}K+\sigma _{16}\right\} p_{j},$ }

\begin{description}
\item {\small \ }

\item[$\left( c\right) $] {\small \ $\overline{k}_{j}=k_{j}+\left\{ \sigma
_{2}w_{2}+\sigma _{3}w_{3}+\sigma _{4}w_{4}\right\} l_{j}$ }
\end{description}

\item {\small $\qquad \qquad +\left\{ -\sigma _{5}w_{2}+\sigma
_{6}w_{3}+\sigma _{7}w_{4}-\sigma _{5}K^{\prime }+\sigma _{6}I^{\prime
}+\sigma _{7}J^{\prime }+\sigma _{19}\right\} m_{j}$ }

\item {\small $\qquad \qquad +\left\{ \sigma _{6}w_{2}-\sigma
_{8}w_{3}-\sigma _{9}w_{4}+\sigma _{6}K^{\prime }-\sigma _{8}I^{\prime
}-\sigma _{9}J^{\prime }-\sigma _{21}\right\} n_{j}$ }

\item {\small $\qquad \qquad +\left\{ \sigma _{7}w_{2}-\sigma
_{9}w_{3}-\sigma _{10}w_{4}+\sigma _{7}K^{\prime }-\sigma _{9}I^{\prime
}-\sigma _{10}J^{\prime }+\sigma _{19}\right\} p_{j}.$\qquad \qquad }
\end{description}

{\small Thus the adapted components $\overline{h}_{\alpha },\overline{j}%
_{\alpha }, \overline{k}_{\alpha }$ of $\overline{h}_{i},\overline{j}_{i},%
\overline{k} _{i}$ in $\overline{F}^{4}=\left( M^{4},\overline{L}\left(
x,y\right) \right) $ are given by }

\begin{description}
\item {\small \ }

\item[$\left( 3.6\right) \left( a\right) $] {\small $\overline{h}%
_{1}=e^{-\sigma }\left\{ h_{1}+\sigma _{2}u_{2}+\sigma _{3}u_{3}+\sigma
_{4}u_{4}\right\} ,$ }

\item {\small $\qquad \qquad \overline{h}_{2}=e^{-\sigma }\left\{
h_{2}-\sigma _{5}u_{2}+\sigma _{6}u_{3}+\sigma _{7}u_{4}+\sigma _{5}\left(
J+J^{\prime }\right) +\sigma _{6}I+\sigma _{7}K^{\prime }-\sigma
_{12}\right\} ,$ }

\item {\small $\qquad \qquad \overline{h}_{3}=e^{-\sigma }\left\{
h_{3}+\sigma _{6}v_{2}-\sigma _{8}v_{3}-\sigma _{9}v_{4}-\sigma _{6}\left(
J+J^{\prime }\right) -\sigma _{8}I-\sigma _{9}K^{\prime }+\sigma
_{14}\right\} ,$ }

\item {\small $\qquad \qquad \overline{h}_{4}=e^{-\sigma }\left\{
h_{4}+\sigma _{7}u_{2}-\sigma _{9}u_{3}-\sigma _{10}u_{4}-\sigma _{7}\left(
J+J^{\prime }\right) -\sigma _{9}I-\sigma _{10}K^{\prime }+\sigma
_{15}\right\} .$ }

\begin{description}
\item {\small \ }

\item[$\left( b\right) $] {\small $\overline{j}_{1}=e^{-\sigma }\left\{
j_{1}+\sigma _{2}v_{2}+\sigma _{3}v_{3}+\sigma _{4}v_{4}\right\} ,$ }
\end{description}

\item {\small $\qquad \qquad \overline{j}_{2}=e^{-\sigma }\left\{
j_{2}-\sigma _{5}v_{2}+\sigma _{6}v_{3}+\sigma _{7}v_{4}+\sigma _{5}\left(
H^{\prime }+I^{\prime }\right) +\sigma _{6}K^{\prime }+\sigma _{7}K-\sigma
_{13}\right\} ,$ }

\item {\small $\qquad \qquad \overline{j}_{3}=e^{-\sigma }\left\{
j_{3}+\sigma _{6}v_{2}-\sigma _{8}v_{3}-\sigma _{9}v_{4}-\sigma _{6}\left(
H^{\prime }+I^{\prime }\right) -\sigma _{8}K^{\prime }-\sigma _{9}K+\sigma
_{15}\right\} ,$ }

\item {\small $\qquad \qquad \overline{j}_{4}=e^{-\sigma }\left\{
j_{4}+\sigma _{7}v_{2}-\sigma _{9}v_{3}-\sigma _{10}v_{4}-\sigma _{7}\left(
H^{\prime }+I^{\prime }\right) -\sigma _{9}K^{\prime }-\sigma _{10}K+\sigma
_{16}\right\} .$ }

\begin{description}
\item {\small \ }

\item[$\left( c\right) $] {\small $\overline{k}_{1}=e^{-\sigma }\left\{
k_{1}+\sigma _{2}w_{2}+\sigma _{3}w_{3}+\sigma _{4}w_{4}\right\} ,$ }

\item {\small $\qquad \overline{k}_{2}=e^{-\sigma }\left\{ k_{2}-\sigma
_{5}w_{2}+\sigma _{6}w_{3}+\sigma _{7}w_{4}-\sigma _{5}K^{\prime }+\sigma
_{6}I^{\prime }+\sigma _{7}J^{\prime }-\sigma _{19}\right\} ,$ }

\item {\small $\qquad \overline{k}_{3}=e^{-\sigma }\left\{ k_{3}+\sigma
_{6}w_{2}-\sigma _{8}w_{3}-\sigma _{9}w_{4}+\sigma _{6}K^{\prime }-\sigma
_{8}I^{\prime }-\sigma _{9}J^{\prime }-\sigma _{21}\right\} ,$ }

\item {\small $\qquad \overline{k}_{4}=e^{-\sigma }\left\{ k_{4}+\sigma
_{7}w_{2}-\sigma _{9}w_{3}-\sigma _{10}w_{4}+\sigma _{7}K^{\prime }-\sigma
_{9}I^{\prime }-\sigma _{10}J^{\prime }+\sigma _{19}\right\} .$ }
\end{description}
\end{description}

\section{\textbf{Conformally flat Landsberg space}}

{\small It is well known that the Berwald spaces are characterized by $%
C_{ijk\mid h}=0$ and Landsberg spaces are characterized by $C_{ijk\mid 0}=0$ 
}${\small (\left[ 2],[14\right] )}${\small $,$ where the index `0' denotes
the transvection by the supporting element $y^{i}.$ Also, one knows that
every Berwald space is a Landsberg space but the converse is not necessarily
true. In }$\left( {\small \left[ 13],[15],[21\right] }\right) ${\small $,$
it has been shown that a Landsberg space becomes a Berwald space under some
special conditions. A lot of work has been done on the conformal
transformation and conformally flat manifolds in Euclidean, hyperbolic,
Riemannian and Finsler geometry. Here, we mention some of the important
results. In $([10],[11])$, N. H. Kuiper studied conformally flat spaces in
large and also for Euclidean spaces of dimension greater than 2. M. Gromov
et. al $[5]$, R. S. Kulkarni$[12]$ have also worked on conformally flat
manifolds. The notion of conformal transformation in Finsler spaces was
introduced by M. Hashiguchi $[6]$. In particular, in dimension }${\small 4}$%
{\small , some of the main contributions are due to M. F. Atiyah et. al $[3],
$ S. Y. Chang et.al $[4]$, M. Iori et. al }${\small [7]}${\small , M.
Kalafat $[8]$, and S. Akubulut et. al $[1]$. In $[9]$, Michael Kapovich has
discussed conformally flat metrics on Riemannian 4-manifolds, which itself
says that there is a lot of possibility to extend these ideas in Finsler
geometry. Motivated by these results, here we obtain the conditions under
which a four-dimensional Landsberg space becomes a Berwald space. We use the
following definition and theorems to obtain the required conditions. }

\begin{definition}
$\left[ 2\right] ${\small A Finsler space $F^{n}$ is called conformally flat
if $F^{n}$ is conformal to a locally Minkowski space. }
\end{definition}

\begin{theorem}
$ \left[ 20\right] $ {\small A Finsler space $F^{4}$ with non-zero $C$ is a Berwald space if and
only if the h-connection vectors $h_{i},j_{i},k_{i}$ vanish and all the main
scalars are $h$-covariant constant. }
\end{theorem}

\begin{theorem}
$ \left[ 20\right] ${\small A Finsler space $F^{4}$ with non-zero C is a Landsberg space if and
only if the h-connection vectors $h_{i},j_{i},k_{i}$ are orthogonal to the
supporting element $y^{i}$ i.e., $h_{1}=j_{1}=k_{1}=0$ and $%
H_{,1}=I_{,1}=J_{,1}=K_{,1}=H_{,1}^{\prime }=I_{,1}^{\prime }=J_{,1}^{\prime
}=K_{,1}^{\prime }=0.$ }
\end{theorem}

{\small If the four-dimensional Finsler space $\overline{F}^{4}=\left( M,%
\overline{L}\right) $ is conformal to the Finsler space $F^{4}=\left(
M,L\right) ,$ the main scalars $\overline{H},\overline{I},\overline{J},%
\overline{K},\overline{H}^{\prime },\overline{I}^{\prime },\overline{J}%
^{\prime },\overline{K\text{ }}^{\prime }$ of $\overline{F}^{4}$ coincide
with the main scalars $H,I,J,K,H^{\prime },I^{\prime },J^{\prime },K^{\prime
}$ of $F^{4}.$ In particular we must notice that the main scalars $%
H,I,J,K,H^{\prime },I^{\prime },J^{\prime },K^{\prime }$ of $F^{4}$ and
h-connection vectors $h_{i},j_{i},k_{i}$ are functions of the variable y$%
^{i} $ alone. Firstly, we assume that the Finsler space $\overline{F}%
^{4}=\left( M,\overline{L}\right) $ is a Landsberg space. Then from Theorem
3, it follows that }

\begin{description}
\item[$\left( 4.1\right) $] {\small $\overline{H}_{,1}=\overline{I}_{,1}=%
\overline{J}_{,1}=\overline{K}_{,1}=\overline{H}_{,1}^{\prime }=\overline{I}%
_{,1}^{\prime }=\overline{J}_{,1}^{\prime }=\overline{K}_{,1}^{\prime }=0$; $%
\overline{h}_{1}=\overline{j}_{1}=\overline{k}_{1}=0.$ }
\end{description}

{\small The scalar $\overline{H}_{,1}$ can be written in terms of Miron
frame as }

{\small $\overline{H}_{,1}=\overline{H}_{,k}\overline{l}^{k}=(\frac{\partial 
\overline{H}}{\partial x^{k}}-\overline{G}_{k}^{r}\frac{\partial \overline{H}%
}{\partial y^{r}})\overline{l}^{k}=-\overline{G}_{k}^{r}\frac{\partial 
\overline{H}}{\partial y^{r}}\overline{l}^{k}=-\overline{G}_{k}^{r}\frac{%
\partial H}{\partial y^{r}}e^{-\sigma }l^{k}$ }

{\small Making use of (3.4), we get }

{\small $\overline{H}_{,1}=-Le^{-\sigma }(\sigma _{1}H\mid _{r}l^{r}-\sigma
_{2}H\mid _{r}m^{r}-\sigma _{3}H\mid _{r}n^{r}-\sigma _{4}H\mid _{r}p^{r}).$ 
}

{\small Therefore, we have }

\begin{description}
\item[$\left( 4.2\right) $] {\small $\overline{H}_{,1}=(\sigma
_{2}H;_{2}+\sigma _{3}H;_{3}+\sigma _{4}H;_{4})e^{-\sigma },\overline{I}%
_{,1}=(\sigma _{2}I;_{2}+\sigma _{3}I;_{3}+\sigma _{4}I;_{4})e^{-\sigma },$ }

\item {\small $\qquad \overline{J}_{,1}=(\sigma _{2}J;_{2}+\sigma
_{3}J;_{3}+\sigma _{4}J;_{4})e^{-\sigma },\overline{K}_{,1}=(\sigma
_{2}K;_{2}+\sigma _{3}K;_{3}+\sigma _{4}K;_{4})e^{-\sigma },$ }

\item {\small $\qquad \overline{H^{\prime }}_{,1}=(\sigma _{2}H^{\prime
};_{2}+\sigma _{3}H^{\prime };_{3}+\sigma _{4}H^{\prime };_{4})e^{-\sigma },%
\overline{I^{\prime }}_{,1}=(\sigma _{2}I^{\prime };_{2}+\sigma
_{3}I^{\prime };_{3}+\sigma _{4}I^{\prime };_{4})e^{-\sigma },$ }

\item {\small $\qquad \overline{J^{\prime }}_{,1}=(\sigma _{2}J^{\prime
};_{2}+\sigma _{3}J^{\prime };_{3}+\sigma _{4}J^{\prime };_{4})e^{-\sigma },%
\overline{K^{\prime }}_{,1}=(\sigma _{2}K^{\prime };_{2}+\sigma
_{3}K^{\prime };_{3}+\sigma _{4}K^{\prime };_{4})e^{-\sigma }.$ }
\end{description}

{\small Therefore, from (4.1), (4.2) and (3.6), we get }

\begin{description}
\item[$\left( 4.3\right) $] {\small $\sigma _{2}H;_{2}+\sigma
_{3}H;_{3}+\sigma _{4}H;_{4}=0,\qquad \sigma _{2}I;_{2}+\sigma
_{3}I;_{3}+\sigma _{4}I;_{4}=0,$ }

\item {\small $\qquad \sigma _{2}J;_{2}+\sigma _{3}J;_{3}+\sigma
_{4}J;_{4}=0,\qquad \sigma _{2}K;_{2}+\sigma _{3}K;_{3}+\sigma _{4}K;_{4}=0,$
}

\item {\small $\qquad \sigma _{2}H^{\prime };_{2}+\sigma _{3}H^{\prime
};_{3}+\sigma _{4}H^{\prime };_{4}=0,\qquad \sigma _{2}I^{\prime
};_{2}+\sigma _{3}I^{\prime };_{3}+\sigma _{4}I^{\prime };_{4}=0,$ }

\item {\small $\qquad \sigma _{2}J^{\prime };_{2}+\sigma _{3}J^{\prime
};_{3}+\sigma _{4}J^{\prime };_{4}=0,\qquad \sigma _{2}K^{\prime
};_{2}+\sigma _{3}K^{\prime };_{3}+\sigma _{4}K^{\prime };_{4}=0$ }
\end{description}

{\small and }

\begin{description}
\item {\small $\qquad h_{1}+\sigma _{2}u_{2}+\sigma _{3}u_{3}+\sigma
_{4}u_{4}=0,\qquad j_{1}+\sigma _{2}v_{2}+\sigma _{3}v_{3}+\sigma
_{4}v_{4}=0,$ }

\item {\small $\qquad k_{1}+\sigma _{2}w_{2}+\sigma _{3}w_{3}+\sigma
_{4}w_{4}=0.$ }
\end{description}

{\small Now, we prove that $\sigma _{2},\sigma _{3}$ and $\sigma _{4}$ never
vanish simultaneously, for non homothetic transformation. If possible let us
assume that $\sigma _{2}=0,\sigma _{3}=0$ and $\sigma _{4}=0$, then $\sigma
_{i}=\sigma _{_{1}}l_{i}+\sigma _{_{2}}m_{i}+\sigma _{_{3}}n_{i}+\sigma
_{_{4}}p_{i}$ gives $\sigma _{i}=\sigma _{_{1}}l_{i}.$ Differentiating this
with respect to $y^{i}$, we get $0=\left( \dot{\partial _{j}}\sigma
_{1}\right) l_{i}+\sigma _{1}\dot{\partial _{j}}l_{i}=\left( \dot{\partial
_{j}}\sigma _{1}\right) l_{i}+\sigma _{1}l_{i}\mid _{j},$ which in view of $%
\left( 2.1\right)(b)$ gives }

\begin{description}
\item[$\left( 4.4\right) $] {\small $\frac{\sigma _{1}}{L}\left(
m_{i}m_{j}+n_{i}n_{j}+p_{i}p_{j}\right) =-\left( \dot{\partial _{j}}\sigma
_{1}\right) l_{i}$ }
\end{description}

{\small Since L. H. S. of the equation $\left( 4.4\right) $ is symmetric in
i and j, we have $\left( \dot{\partial _{i}}\sigma _{1}\right) l_{j}=\left( 
\dot{\partial _{j}}\sigma _{1}\right) l_{i}.$ Contracting this equation by $%
l^{j}$, we get $\dot{\partial _{i}}\sigma _{1}=\left( \dot{\partial _{j}}%
\sigma _{1}\right) l_{i}l^{j}.$ Since $\sigma _{1}$ is positively
homogeneous of degree zero in y$^{i}$, we get $\left( \dot{\partial _{j}}%
\sigma _{1}\right) l^{j}=0,$ which implies $\dot{\partial _{i}}\sigma _{1}=0.
$ Thus, the equation $\left( 4.4\right) $ shows that $\sigma _{1}=0.$ Hence $%
\sigma _{i}=0$ $(i=1,...,4);$ which shows that $\sigma $ is constant, i.e.,
the transformation is homothetic. Hence, we conclude that, for a
non-homothetic transformation $\sigma _{2},\sigma _{3}$ and $\sigma _{4}$
never vanish simultaneously. So we consider the following cases for
non-homothetic transformation. }

\begin{description}
\item[Case(i)] {\small When $\sigma _{2}\neq 0,\sigma _{3}$ $\neq 0$, $%
\sigma _{4}\neq 0.$ }
\end{description}

{\small In this case, we see that there is no change in equation $\left(
4.3\right) .$ From this, it follows that if the space is a Landsberg space
then (4.3) holds. }

{\small Conversely, if $\left( 4.3\right) $ holds then from (4.2) and (3.6),
we get (4.1). So, $\left( M,\overline{L}\right) $ is a Landsberg space. }

\begin{description}
\item[Case(ii)] {\small When $\sigma _{2}\neq 0,\sigma _{3}$ $\neq 0$, $%
\sigma _{4}=0.$ }
\end{description}

{\small In this case, from the equation (4.3), we get }

\begin{description}
\item[$\left( 4.5\right) $] {\small $\frac{H;_{2}}{H;_{3}}=\frac{I;_{2}}{%
I;_{3}}=\frac{J;_{2}}{J;_{3}}=\frac{K;_{2}}{K;_{3}}=\frac{H^{\prime };_{2}}{%
H^{\prime };_{3}}=\frac{I^{\prime };_{2}}{I^{\prime };_{3}}=\frac{J^{\prime
};_{2}}{J^{\prime };_{3}}=\frac{K^{\prime };_{2}}{K^{\prime };_{3}}=-\frac{%
\sigma _{3}}{\sigma _{2}}$ }
\end{description}

{\small and }

\begin{description}
\item {\small $\qquad h_{1}=-\left( \sigma _{2}u_{2}+\sigma _{3}u_{3}\right)
,j_{1}=-\left( \sigma _{2}v_{2}+\sigma _{3}v_{3}\right) ,k_{1}=-(\sigma
_{2}w_{2}+\sigma _{3}w_{3}).$ }
\end{description}

{\small Conversely, if (4.5) holds then from (3.6), (4.2) and (4.3), we get
(4.1). So, (M, $\overline{L})$ is a Landsberg space. }

\begin{description}
\item[Case(iii)] {\small When $\sigma _{2}\neq 0,\sigma _{3}$ $=0$, $\sigma
_{4}\neq 0.$ }
\end{description}

{\small In this case, from the equation (4.3), we get }

\begin{description}
\item[$\left( 4.6\right) $] {\small $\frac{H;_{2}}{H;_{4}}=\frac{I;_{2}}{%
I;_{4}}=\frac{J;_{2}}{J;4}=\frac{K;_{2}}{K;_{4}}=\frac{H^{\prime };_{2}}{%
H^{\prime };_{4}}=\frac{I^{\prime };_{2}}{I^{\prime };_{4}}=\frac{J^{\prime
};_{2}}{J^{\prime };_{4}}=\frac{K^{\prime };_{2}}{K^{\prime };_{4}}=-\frac{%
\sigma _{4}}{\sigma _{2}}$ }
\end{description}

{\small and }

\begin{description}
\item {\small $\qquad h_{1}=-\left( \sigma _{2}u_{2}+\sigma _{4}u_{4}\right)
,j_{1}=-\left( \sigma _{2}v_{2}+\sigma _{4}v_{4}\right) ,k_{1}=-(\sigma
_{2}w_{2}+\sigma _{4}w_{4}).$ }
\end{description}

{\small Conversely, if (4.6) holds then from (3.6), (4.2) and (4.3), we get
(4.1). So, (M, $\overline{L})$ is a Landsberg space. }

\begin{description}
\item[Case(iv)] {\small When $\sigma _{2}=0,\sigma _{3}$ $\neq 0$, $\sigma
_{4}\neq 0.$ }
\end{description}

{\small In this case, from the equation (4.3), we get }

\begin{description}
\item[$\left( 4.7\right) $] {\small $\frac{H;_{3}}{H;_{4}}=\frac{I;_{3}}{%
I;_{4}}=\frac{J;_{3}}{J;4}=\frac{K;_{3}}{K;_{4}}=\frac{H^{\prime };_{3}}{%
H^{\prime };_{4}}=\frac{I^{\prime };_{3}}{I^{\prime };_{4}}=\frac{J^{\prime
};_{3}}{J^{\prime };_{4}}=\frac{K^{\prime };_{3}}{K^{\prime };_{4}}=-\frac{%
\sigma _{4}}{\sigma _{3}}$ }
\end{description}

{\small and }

\begin{description}
\item {\small $\qquad h_{1}=-\left( \sigma _{3}u_{3}+\sigma _{4}u_{4}\right)
,j_{1}=-\left( \sigma _{3}v_{3}+\sigma _{4}v_{4}\right) ,k_{1}=-(\sigma
_{3}w_{3}+\sigma _{4}w_{4}).$ }
\end{description}

{\small Conversely, if (4.7) holds then from (3.6), (4.2) and (4.3), we get
(4.1). So, $(M,\overline{L})$ is a Landsberg space.}

\begin{description}
\item[Case(v)] {\small When $\sigma _{2}\neq 0,\sigma _{3}$=$0$, $\sigma
_{4}=0.$ }
\end{description}

{\small In this case, from the equation (4.3), we get }

\begin{description}
\item[$\left( 4.8\right) $] {\small $H;_{2}=I;_{2}=J;_{2}=K;_{2}=H^{\prime
};_{2}=I^{\prime };_{2}=J^{\prime };_{2}=K^{\prime };_{2}=0$ }
\end{description}

{\small and }

\begin{description}
\item {\small $\qquad h_{1}=-\sigma _{2}u_{2},\qquad j_{1}=-\sigma
_{2}v_{2},\qquad k_{1}=-\sigma _{2}w_{2}.$ }
\end{description}

{\small Conversely, if (4.8) holds then from (3.6), (4.2) and (4.3), we get
(4.1). So, (M, $\overline{L})$ is a Landsberg space. }

\begin{description}
\item[Case(vi)] {\small When $\sigma _{2}=0,\sigma _{3}\neq 0$, $\sigma
_{4}=0.$ }
\end{description}

{\small In this case, from the equation (4.3), we get }

\begin{description}
\item[$\left( 4.9\right) $] {\small $H;_{3}=I;_{3}=J;_{3}=K;_{3}=H^{\prime
};_{3}=I^{\prime };_{3}=J^{\prime };_{3}=K^{\prime };_{3}=0$ }
\end{description}

{\small and }

\begin{description}
\item {\small $\qquad h_{1}=-\sigma _{3}u_{3},\qquad j_{1}=-\sigma
_{3}v_{3},\qquad k_{1}=-\sigma _{3}w_{3}.$ }
\end{description}

{\small Conversely, if (4.9) holds then from (3.6), (4.2) and (4.3), we get
(4.1). So, (M, $\overline{L})$ is a Landsberg space. }

\begin{description}
\item[Case(vii)] {\small When $\sigma _{2}=0,\sigma _{3}=0$, $\sigma
_{4}\neq 0.$ }
\end{description}

{\small In this case, from the equation (4.3), we get }

\begin{description}
\item[$\left( 4.10\right) $] {\small $H;_{4}=I;_{4}=J;_{4}=K;_{4}=H^{\prime
};_{4}=I^{\prime };_{4}=J^{\prime };_{4}=K^{\prime };_{4}=0$ }
\end{description}

{\small and }

\begin{description}
\item {\small $\qquad h_{1}=-\sigma _{4}u_{4},\qquad j_{1}=-\sigma
_{4}v_{4},\qquad k_{1}=-\sigma _{4}w_{4}.$ }
\end{description}

{\small Conversely, if (4.10) holds then from (3.6), (4.2) and (4.3), we get
(4.1). So, (M, $\overline{L})$ is a Landsberg space. Hence we have the
following: }

\begin{theorem}
{\small A four-dimensional Landsberg space is $\sigma -$ conformally flat if
and only if one of the following conditions is satisfied:}

\begin{description}
\item 1) the equation (4.3) holds for {\small $\sigma _{2}\neq 0,\sigma _{3}$
$\neq 0$, $\sigma _{4}\neq 0.$}

\item 2) {\small $\sigma _{i}$ is orthogonal to $p^{i}$ and the equation
(4.5) holds for $\sigma _{2}\neq 0,\sigma _{3}$ $\neq 0$, $\sigma _{4}=0.$}

\item 3) {\small $\sigma _{i}$ is orthogonal to $n^{i}$ and the equation
(4.6) holds for $\sigma _{2}\neq 0,\sigma _{3}$ $=0$, $\sigma _{4}\neq 0.$}

\item 4) {\small $\sigma _{i}$ is orthogonal to $m^{i}$ and the equation
(4.7) holds for $\sigma _{2}=0,\sigma _{3}$ $\neq 0$, $\sigma _{4}\neq 0.$}

\item 5) {\small $\sigma _{i}$ is orthogonal to $n^{i}$ and $p^{i}$ and the
equation (4.8) holds for $\sigma _{2}\neq 0,\sigma _{3}$=$0$, $\sigma
_{4}=0. $}

\item 6) {\small $\sigma _{i}$ is orthogonal to $m^{i}$ and $p^{i}$ and the
equation (4.9) holds for $\sigma _{2}=0,\sigma _{3}\neq 0$, $\sigma _{4}=0.$}

\item 7) {\small $\sigma _{i}$ is orthogonal to $m^{i}$ and $n^{i},$ and the
equation (4.10) holds for $\sigma _{2}=0,\sigma _{3}=0$, $\sigma _{4}\neq 0.$
}
\end{description}
\end{theorem}

\section{\textbf{Conformally flat Berwald space}}

{\small A Finsler metric is called conformally Berwald Finsler metric if it
is conformally related to a Berwald metric. In particular, if a Finsler
metric is conformally related to a Minkowski metric, then it is called
conformally flat Finsler metric. More precisely, a Finsler manifold $\left(
M,L\right) $ is called conformally flat or $\sigma $- conformally flat, if
there exists an atlas whose changes of coordinates are conformal
diffeomorphism to open sets in some Minkowski space. We consider the case
when the Finsler space $\overline{F}^{4}=\left( M,\overline{L}\right) $ is a
Berwald space. We shall write $\overline{H}_{\bot k}=(\frac{\partial 
\overline{H}}{\partial x^{k}}-\overline{G}_{k}^{r}\frac{\partial \overline{H}%
}{\partial y^{r}}).$ Since the main scalars $H,I,J,K,H^{\prime },I^{\prime
},J^{\prime },K^{\prime }$ of $F^{4}$ and h-connection vectors $%
h_{i},j_{i},k_{i}$ are functions of the variable y$^{i}$ alone, the above
equation is equivalent to $\overline{H}_{\bot k}=-\overline{G}_{k}^{r}\frac{%
\partial H}{\partial y^{r}}.$ Making use of (3.4), we get }

\begin{description}
\item[$\left( 5.1\right) $] {\small $\overline{H}_{\perp k}=-L\{l^{r}\left(
\sigma _{1}l_{k}+\sigma _{2}m_{k}+\sigma _{3}n_{k}+\sigma _{4}p_{k}\right)
-m^{r}\left( \sigma _{2}l_{k}-\sigma _{5}m_{k}+\sigma _{6}n_{k}+\sigma
_{7}p_{k}\right) $ }

\item {\small $\qquad -n^{r}\left( \sigma _{3}l_{k}+\sigma _{6}m_{k}-\sigma
_{8}n_{k}-\sigma _{9}p_{k}\right) -p^{r}\left( \sigma _{4}l_{k}+\sigma
_{7}m_{k}-\sigma _{9}n_{k}-\sigma _{10}p_{k}\right) \}H\mid _{r}.$ }
\end{description}

{\small Since $H\mid
_{r}=L^{-1}(H;_{1}l_{r}+H;_{2}m_{r}+H;_{3}n_{r}+H;_{4}p_{r})$ and $H;_{1}=0,$
we have }

\begin{description}
\item[$\left( 5.2\right) $] {\small $\overline{H}_{\perp k}=H;_{2}\left(
\sigma _{2}l_{k}-\sigma _{5}m_{k}+\sigma _{6}n_{k}+\sigma _{7}p_{k}\right)
+H;_{3}\left( \sigma _{3}l_{k}+\sigma _{6}m_{k}-\sigma _{8}n_{k}-\sigma
_{9}p_{k}\right) $ }

\item {\small $\qquad \qquad +H;_{4}\left( \sigma _{6}l_{k}+\sigma
_{7}m_{k}-\sigma _{9}n_{k}-\sigma _{10}p_{k}\right) ,$ }

\item {\small $\qquad \overline{I}_{\perp k}=I;_{2}\left( \sigma
_{2}l_{k}-\sigma _{5}m_{k}+\sigma _{6}n_{k}+\sigma _{7}p_{k}\right)
+I;_{3}\left( \sigma _{3}l_{k}+\sigma _{6}m_{k}-\sigma _{8}n_{k}-\sigma
_{9}p_{k}\right) $ }

\item {\small $\qquad \qquad +I;_{4}\left( \sigma _{6}l_{k}+\sigma
_{7}m_{k}-\sigma _{9}n_{k}-\sigma _{10}p_{k}\right) ,$ }

\item {\small $\qquad \overline{J}_{\perp k}=J;_{2}\left( \sigma
_{2}l_{k}-\sigma _{5}m_{k}+\sigma _{6}n_{k}+\sigma _{7}p_{k}\right)
+J;_{3}\left( \sigma _{3}l_{k}+\sigma _{6}m_{k}-\sigma _{8}n_{k}-\sigma
_{9}p_{k}\right) $ }

\item {\small $\qquad \qquad +J;_{4}\left( \sigma _{6}l_{k}+\sigma
_{7}m_{k}-\sigma _{9}n_{k}-\sigma _{10}p_{k}\right) ,$ }

\item {\small $\qquad \overline{K}_{\perp k}=K;_{2}\left( \sigma
_{2}l_{k}-\sigma _{5}m_{k}+\sigma _{6}n_{k}+\sigma _{7}p_{k}\right)
+K;_{3}\left( \sigma _{3}l_{k}+\sigma _{6}m_{k}-\sigma _{8}n_{k}-\sigma
_{9}p_{k}\right) $ }

\item {\small $\qquad \qquad +K;_{4}\left( \sigma _{6}l_{k}+\sigma
_{7}m_{k}-\sigma _{9}n_{k}-\sigma _{10}p_{k}\right) ,$ }

\item {\small $\qquad \overline{H^{\prime }}_{\perp k}=H^{\prime
};_{2}\left( \sigma _{2}l_{k}-\sigma _{5}m_{k}+\sigma _{6}n_{k}+\sigma
_{7}p_{k}\right) +H^{\prime };_{3}\left( \sigma _{3}l_{k}+\sigma
_{6}m_{k}-\sigma _{8}n_{k}-\sigma _{9}p_{k}\right) $ }

\item {\small $\qquad \qquad +H^{\prime };_{4}\left( \sigma _{6}l_{k}+\sigma
_{7}m_{k}-\sigma _{9}n_{k}-\sigma _{10}p_{k}\right) ,$ }

\item {\small $\qquad \overline{I^{\prime }}_{\perp k}=I^{\prime
};_{2}\left( \sigma _{2}l_{k}-\sigma _{5}m_{k}+\sigma _{6}n_{k}+\sigma
_{7}p_{k}\right) +I^{\prime };_{3}\left( \sigma _{3}l_{k}+\sigma
_{6}m_{k}-\sigma _{8}n_{k}-\sigma _{9}p_{k}\right) $ }

\item {\small $\qquad \qquad +I^{\prime };_{4}\left( \sigma _{6}l_{k}+\sigma
_{7}m_{k}-\sigma _{9}n_{k}-\sigma _{10}p_{k}\right) ,$ }

\item {\small $\qquad \overline{J^{\prime }}_{\perp k}=J^{\prime
};_{2}\left( \sigma _{2}l_{k}-\sigma _{5}m_{k}+\sigma _{6}n_{k}+\sigma
_{7}p_{k}\right) +J^{\prime };_{3}\left( \sigma _{3}l_{k}+\sigma
_{6}m_{k}-\sigma _{8}n_{k}-\sigma _{9}p_{k}\right) $ }

\item {\small $\qquad \qquad +J^{\prime };_{4}\left( \sigma _{6}l_{k}+\sigma
_{7}m_{k}-\sigma _{9}n_{k}-\sigma _{10}p_{k}\right) ,$ }

\item {\small $\qquad \overline{K^{\prime }}_{\perp k}=K^{\prime
};_{2}\left( \sigma _{2}l_{k}-\sigma _{5}m_{k}+\sigma _{6}n_{k}+\sigma
_{7}p_{k}\right) +K^{\prime };_{3}\left( \sigma _{3}l_{k}+\sigma
_{6}m_{k}-\sigma _{8}n_{k}-\sigma _{9}p_{k}\right) $ }

\item {\small $\qquad \qquad +K^{\prime };_{4}\left( \sigma _{6}l_{k}+\sigma
_{7}m_{k}-\sigma _{9}n_{k}-\sigma _{10}p_{k}\right) .$ }
\end{description}

{\small Now, we discuss all the seven cases which have been discussed in the
previous section. }

\begin{description}
\item[Case(i)] {\small When $\sigma _{2}\neq 0,\sigma _{3}\neq 0,\sigma
_{4}\neq 0.$ }
\end{description}

{\small If the four dimensional Landsberg space $\left( M,\overline{L}%
\right) $ is conformally flat, then from the equations (5.2), (3.5) and
(4.3), we get\qquad }

\begin{description}
\item[$\left( 5.3\right) $] {\small $\overline{H}_{\perp k}=\left( -\sigma
_{5}H;_{2}+\sigma _{6}H;_{3}+\sigma _{7}H;_{4}\right) m_{k}+\left( \sigma
_{6}H;_{2}-\sigma _{8}H;_{3}-\sigma _{9}H;_{4}\right) n_{k}$ }

\item {\small $\qquad \qquad +\left( \sigma _{7}H;_{2}-\sigma
_{9}H;_{3}-\sigma _{10}H;_{4}\right) p_{k},$ }

\item {\small $\overline{I}_{\perp k}=\left( -\sigma _{5}I;_{2}+\sigma
_{6}I;_{3}+\sigma _{7}I;_{4}\right) m_{k}+\left( \sigma _{6}I;_{2}-\sigma
_{8}I;_{3}-\sigma _{9}I;_{4}\right) n_{k}$ }

\item {\small $\qquad \qquad +\left( \sigma _{7}I;_{2}-\sigma
_{9}I;_{3}-\sigma _{10}I;_{4}\right) p_{k},$ }

\item {\small $\overline{J}_{\perp k}=\left( -\sigma _{5}J;_{2}+\sigma
_{6}J;_{3}+\sigma _{7}J;_{4}\right) m_{k}+\left( \sigma _{6}J;_{2}-\sigma
_{8}J;_{3}-\sigma _{9}J;_{4}\right) n_{k}$ }

\item {\small $\qquad \qquad +\left( \sigma _{7}J;_{2}-\sigma
_{9}J;_{3}-\sigma _{10}J;_{4}\right) p_{k},$ }

\item {\small $\overline{K}_{\perp k}=\left( -\sigma _{5}K;_{2}+\sigma
_{6}K;_{3}+\sigma _{7}K;_{4}\right) m_{k}+\left( \sigma _{6}K;_{2}-\sigma
_{8}K;_{3}-\sigma _{9}K;_{4}\right) n_{k}$ }

\item {\small $\qquad \qquad +\left( \sigma _{7}K;_{2}-\sigma
_{9}K;_{3}-\sigma _{10}K;_{4}\right) p_{k},$ }

\item {\small $\overline{H^{\prime }}_{\perp k}=\left( -\sigma _{5}H^{\prime
};_{2}+\sigma _{6}H^{\prime };_{3}+\sigma _{7}H^{\prime };_{4}\right)
m_{k}+\left( \sigma _{6}H^{\prime };_{2}-\sigma _{8}H^{\prime };_{3}-\sigma
_{9}H^{\prime };_{4}\right) n_{k}$ }

\item {\small $\qquad \qquad +\left( \sigma _{7}H^{\prime };_{2}-\sigma
_{9}H^{\prime };_{3}-\sigma _{10}H^{\prime };_{4}\right) p_{k},$ }

\item {\small $\overline{I^{\prime }}_{\perp k}=\left( -\sigma _{5}I^{\prime
};_{2}+\sigma _{6}I^{\prime };_{3}+\sigma _{7}I^{\prime };_{4}\right)
m_{k}+\left( \sigma _{6}I^{\prime };_{2}-\sigma _{8}I^{\prime };_{3}-\sigma
_{9}I^{\prime };_{4}\right) n_{k}$ }

\item {\small $\qquad \qquad +\left( \sigma _{7}I^{\prime };_{2}-\sigma
_{9}I^{\prime };_{3}-\sigma _{10}I^{\prime };_{4}\right) p_{k},$ }

\item {\small $\overline{J^{\prime }}_{\perp k}=\left( -\sigma _{5}J^{\prime
};_{2}+\sigma _{6}J^{\prime };_{3}+\sigma _{7}J^{\prime };_{4}\right)
m_{k}+\left( \sigma _{6}J^{\prime };_{2}-\sigma _{8}J^{\prime };_{3}-\sigma
_{9}J^{\prime };_{4}\right) n_{k}$ }

\item {\small $\qquad \qquad +\left( \sigma _{7}J^{\prime };_{2}-\sigma
_{9}J^{\prime };_{3}-\sigma _{10}K^{\prime };_{4}\right) p_{k},$ }

\item {\small $\overline{K^{\prime }}_{\perp k}=\left( -\sigma _{5}K^{\prime
};_{2}+\sigma _{6}K^{\prime };_{3}+\sigma _{7}K^{\prime };_{4}\right)
m_{k}+\left( \sigma _{6}K^{\prime };_{2}-\sigma _{8}K^{\prime };_{3}-\sigma
_{9}K^{\prime };_{4}\right) n_{k}$ }

\item {\small $\qquad \qquad +\left( \sigma _{7}K^{\prime };_{2}-\sigma
_{9}K^{\prime };_{3}-\sigma _{10}K^{\prime };_{4}\right) p_{k}.$ }
\end{description}

{\small and }

\begin{description}
\item {\small $\overline{h}_{j}=[h_{2}+\{-\sigma _{5}u_{2}+\sigma
_{6}u_{3}+\sigma _{7}u_{4}+\sigma _{5}\left( J+J^{\prime }\right) +\sigma
_{6}I+\sigma _{7}K^{\prime }-\sigma _{12}\}]m_{j}$ }

\item {\small $\qquad +[h_{3}+\{\sigma _{6}u_{2}-\sigma _{8}u_{3}-\sigma
_{9}u_{4}-\sigma _{6}\left( J+J^{\prime }\right) -\sigma _{8}I-\sigma
_{9}K^{\prime }+\sigma _{14}\}]n_{j}$ }

\item {\small $\qquad +[h_{4}+\{\sigma _{7}u_{2}-\sigma _{9}u_{3}-\sigma
_{10}u_{4}-\sigma _{7}\left( J+J^{\prime }\right) -\sigma _{9}I-\sigma
_{10}K^{\prime }+\sigma _{15}\}]p_{j},$ }

\item {\small $\overline{j}_{j}=[j_{2}+\{-\sigma _{5}v_{2}+\sigma
_{6}v_{3}+\sigma _{7}v_{4}+\sigma _{5}\left( H^{\prime }+I^{\prime }\right)
+\sigma _{6}K^{\prime }+\sigma _{7}K-\sigma _{13}\}]m_{j}$ }

\item {\small $\qquad +[j_{3}+\{\sigma _{6}v_{2}-\sigma _{8}v_{3}-\sigma
_{9}v_{4}-\sigma _{6}\left( H^{\prime }+I^{\prime }\right) -\sigma
_{8}K^{\prime }-\sigma _{9}K+\sigma _{15}\}]n_{j}$ }

\item {\small $\qquad +[j_{4}+\{\sigma _{7}v_{2}-\sigma _{9}v_{3}-\sigma
_{10}v_{4}-\sigma _{7}\left( H^{\prime }+I^{\prime }\right) -\sigma
_{9}K^{\prime }-\sigma _{10}K+\sigma _{16}\}]p_{j},$ }

\item {\small \ $\overline{k}_{j}=[k_{2}+\left\{ -\sigma _{5}w_{2}+\sigma
_{6}w_{3}+\sigma _{7}w_{4}-\sigma _{5}K^{\prime }+\sigma _{6}I^{\prime
}+\sigma _{7}J^{\prime }+\sigma _{19}\right\} ]m_{j}$ }

\item {\small $\qquad +[k_{3}+\left\{ \sigma _{6}w_{2}-\sigma
_{8}w_{3}-\sigma _{9}w_{4}+\sigma _{6}K^{\prime }-\sigma _{8}I^{\prime
}-\sigma _{9}J^{\prime }-\sigma _{21}\right\} ]n_{j}$ }

\item {\small $\qquad +[k_{4}+\left\{ \sigma _{7}w_{2}-\sigma
_{9}w_{3}-\sigma _{10}w_{4}+\sigma _{7}K^{\prime }-\sigma _{9}I^{\prime
}-\sigma _{10}J^{\prime }+\sigma _{19}\right\} ]p_{j}.$ }
\end{description}

{\small From Theorem 2, it follows that the space F$^{4}=\left( M,\overline{L%
}\right) $ is a Berwald space if $\overline{H}_{\perp k}=\overline{I}_{\perp
k}=\overline{J}_{\perp k}=\overline{K}_{\perp k}=\overline{H}_{\perp
k}^{\prime }=\overline{I^{\prime }}_{\perp k}=\overline{J^{\prime }}_{\perp
k}=\overline{K^{\prime }}_{\perp k}=0$ and $\overline{h_{i}}=\overline{j_{i}}%
=\overline{k_{i}}=0.$ Therefore from (5.3), it follows that the space $%
\left( M,\overline{L}\right) $ is a Berwald space if }

\begin{description}
\item[$\left( 5.4\right) $] {\small $-\sigma _{5}H;_{2}+\sigma
_{6}H;_{3}+\sigma _{7}H;_{4}=0,\qquad \sigma _{6}H;_{2}-\sigma
_{8}H;_{3}-\sigma _{9}H;_{4}=0,\qquad $ }

\item {\small $\qquad \sigma _{7}H;_{2}-\sigma _{9}H;_{3}-\sigma
_{10}H;_{4}=0,\qquad -\sigma _{5}I;_{2}+\sigma _{6}I;_{3}+\sigma
_{7}I;_{4}=0,$ }

\item {\small $\qquad \sigma _{6}I;_{2}-\sigma _{8}I;_{3}-\sigma
_{9}I;_{4}=0,\qquad \sigma _{7}I;_{2}-\sigma _{9}I;_{3}-\sigma
_{10}I;_{4}=0, $ }

\item {\small $\qquad -\sigma _{5}J;_{2}+\sigma _{6}J;_{3}+\sigma
_{7}J;_{4}=0,\qquad \sigma _{6}J;_{2}-\sigma _{8}J;_{3}-\sigma _{9}J;_{4}=0,$
}

\item {\small $\qquad \sigma _{7}J;_{2}-\sigma _{9}J;_{3}-\sigma
_{10}J;_{4}=0,\qquad -\sigma _{5}K;_{2}+\sigma _{6}K;_{3}+\sigma
_{7}K;_{4}=0,$ }

\item {\small $\qquad \sigma _{6}K;_{2}-\sigma _{8}K;_{3}-\sigma
_{9}K;_{4}=0,\qquad \sigma _{7}K;_{2}-\sigma _{9}K;_{3}-\sigma
_{10}K;_{4}=0, $ }

\item {\small $\qquad -\sigma _{5}H^{\prime };_{2}+\sigma _{6}H^{\prime
};_{3}+\sigma _{7}H^{\prime };_{4}=0,\qquad \sigma _{6}H^{\prime
};_{2}-\sigma _{8}H^{\prime };_{3}-\sigma _{9}H^{\prime };_{4}=0,$ }

\item {\small $\qquad \sigma _{7}H^{\prime };_{2}-\sigma _{9}H^{\prime
};_{3}-\sigma _{10}H^{\prime };_{4}=0,\qquad -\sigma _{5}I^{\prime
};_{2}+\sigma _{6}I^{\prime };_{3}+\sigma _{7}I^{\prime };_{4}=0,$ }

\item {\small $\qquad \sigma _{6}I^{\prime };_{2}-\sigma _{8}I^{\prime
};_{3}-\sigma _{9}I^{\prime };_{4}=0,\qquad \sigma _{7}I^{\prime
};_{2}-\sigma _{9}I^{\prime };_{3}-\sigma _{10}I^{\prime };_{4}=0,$ }

\item {\small $\qquad -\sigma _{5}J^{\prime };_{2}+\sigma _{6}J^{\prime
};_{3}+\sigma _{7}J^{\prime };_{4}=0,\qquad \sigma _{6}J^{\prime
};_{2}-\sigma _{8}J^{\prime };_{3}-\sigma _{9}J^{\prime };_{4}=0,$ }

\item {\small $\qquad \sigma _{7}J^{\prime };_{2}-\sigma _{9}J^{\prime
};_{3}-\sigma _{10}K^{\prime };_{4}=0,\qquad -\sigma _{5}K^{\prime
};_{2}+\sigma _{6}K^{\prime };_{3}+\sigma _{7}K^{\prime };_{4}=0,$ }

\item {\small $\qquad \sigma _{6}K^{\prime };_{2}-\sigma _{8}K^{\prime
};_{3}-\sigma _{9}K^{\prime };_{4}=0,\qquad \sigma _{7}K^{\prime
};_{2}-\sigma _{9}K^{\prime };_{3}-\sigma _{10}K^{\prime };_{4}=0.$ }
\end{description}

{\small and }

\begin{description}
\item {\small $\qquad h_{2}+\{-\sigma _{5}u_{2}+\sigma _{6}u_{3}+\sigma
_{7}u_{4}+\sigma _{5}\left( J+J^{\prime }\right) +\sigma _{6}I+\sigma
_{7}K^{\prime }-\sigma _{12}\}=0,$ }

\item {\small $\qquad h_{3}+\{\sigma _{6}u_{2}-\sigma _{8}u_{3}-\sigma
_{9}u_{4}-\sigma _{6}\left( J+J^{\prime }\right) -\sigma _{8}I-\sigma
_{9}K^{\prime }+\sigma _{14}\}=0,$ }

\item {\small $\qquad h_{4}+\{\sigma _{7}u_{2}-\sigma _{9}u_{3}-\sigma
_{10}u_{4}-\sigma _{7}\left( J+J^{\prime }\right) -\sigma _{9}I-\sigma
_{10}K^{\prime }+\sigma _{15}\}=0,$ }

\item {\small $\qquad j_{2}+\{-\sigma _{5}v_{2}+\sigma _{6}v_{3}+\sigma
_{7}v_{4}+\sigma _{5}\left( H^{\prime }+I^{\prime }\right) +\sigma
_{6}K^{\prime }+\sigma _{7}K-\sigma _{13}\}=0,$ }

\item {\small $\qquad j_{3}+\{\sigma _{6}v_{2}-\sigma _{8}v_{3}-\sigma
_{9}v_{4}-\sigma _{6}\left( H^{\prime }+I^{\prime }\right) -\sigma
_{8}K^{\prime }-\sigma _{9}K+\sigma _{15}\}=0,$ }

\item {\small $\qquad j_{4}+\{\sigma _{7}v_{2}-\sigma _{9}v_{3}-\sigma
_{10}v_{4}-\sigma _{7}\left( H^{\prime }+I^{\prime }\right) -\sigma
_{9}K^{\prime }-\sigma _{10}K+\sigma _{16}\}=0,$ }

\item {\small \qquad\ $k_{2}+\left\{ -\sigma _{5}w_{2}+\sigma
_{6}w_{3}+\sigma _{7}w_{4}-\sigma _{5}K^{\prime }+\sigma _{6}I^{\prime
}+\sigma _{7}J^{\prime }+\sigma _{19}\right\} =0,$ }

\item {\small $\qquad k_{3}+\left\{ \sigma _{6}w_{2}-\sigma _{8}w_{3}-\sigma
_{9}w_{4}+\sigma _{6}K^{\prime }-\sigma _{8}I^{\prime }-\sigma _{9}J^{\prime
}-\sigma _{21}\right\} =0,$ }

\item {\small $\qquad k_{4}+\left\{ \sigma _{7}w_{2}-\sigma _{9}w_{3}-\sigma
_{10}w_{4}+\sigma _{7}K^{\prime }-\sigma _{9}I^{\prime }-\sigma
_{10}J^{\prime }+\sigma _{19}\right\} =0.$ }
\end{description}

{\small Conversely, if (5.4) holds, then from (5.3), we get $\overline{H}%
_{\perp k}=\overline{I}_{\perp k}=\overline{J}_{\perp k}=\overline{K}_{\perp
k}=\overline{H}_{\perp k}^{\prime }=\overline{I^{\prime }}_{\perp k}=%
\overline{J^{\prime }}_{\perp k}=\overline{K^{\prime }}_{\perp k}=0$ and $%
\overline{h_{i}}=\overline{j_{i}}=\overline{k_{i}}=0.$ Hence the space $%
\left( M,\overline{L}\right) $ is a Berwald space. }In the cases\newline
\textbf{(ii)} {\small When $\sigma _{2}\neq 0,\sigma _{3}$ $\neq 0$, $\sigma
_{4}=0,$ }\newline
\textbf{(iii)} {\small When $\sigma _{2}\neq 0,\sigma _{3}$ $=0$, $\sigma
_{4}\neq 0,$ } and\newline
\textbf{(iv)} {\small When $\sigma _{2}=0,\sigma _{3}\neq 0$, $\sigma
_{4}\neq 0,$ }\newline
{\small we see that if a four- dimensional Landsberg space is $\sigma -$%
conformally flat, then from the equations (5.2), (3.5) and (4.5), there is
no change in the equation (5.4).}

\begin{description}
\item[Case(v)] {\small When $\sigma _{2}\neq 0,\sigma _{3}$ $=0$, $\sigma
_{4}=0.$ }
\end{description}

{\small In this case, if a four- dimensional Landsberg space is $\sigma -$%
conformally flat, then from the equations (5.2), (3.5) and (4.8), we see
that (5.4) reduces to }

\begin{description}
\item[$\left( 5.5\right) $] {\small $\frac{H;_{3}}{H;_{4}}=\frac{I;_{3}}{%
I;_{4}}=\frac{J;_{3}}{J;_{4}}=\frac{K;_{3}}{K;_{4}}=\frac{H^{\prime };_{3}}{%
H^{\prime };_{4}}=\frac{I^{\prime };_{3}}{I^{\prime };_{4}}=\frac{J^{\prime
};_{3}}{J^{\prime };_{4}}=\frac{K^{\prime };_{3}}{K^{\prime };_{4}}=-\frac{%
\sigma _{7}}{\sigma _{6}}=-\frac{\sigma _{9}}{\sigma _{8}}=-\frac{\sigma
_{10}}{\sigma _{9}}$ }
\end{description}

{\small and }

\begin{description}
\item {\small $\qquad h_{2}+\{-\sigma _{5}u_{2}+\sigma _{6}u_{3}+\sigma
_{7}u_{4}+\sigma _{5}\left( J+J^{\prime }\right) +\sigma _{6}I+\sigma
_{7}K^{\prime }-\sigma _{12}\}=0,$ }

\item {\small $\qquad h_{3}+\{\sigma _{6}u_{2}-\sigma _{8}u_{3}-\sigma
_{9}u_{4}-\sigma _{6}\left( J+J^{\prime }\right) -\sigma _{8}I-\sigma
_{9}K^{\prime }+\sigma _{14}\}=0,$ }

\item {\small $\qquad h_{4}+\{\sigma _{7}u_{2}-\sigma _{9}u_{3}-\sigma
_{10}u_{4}-\sigma _{7}\left( J+J^{\prime }\right) -\sigma _{9}I-\sigma
_{10}K^{\prime }+\sigma _{15}\}=0,$ }

\item {\small $\qquad j_{2}+\{-\sigma _{5}v_{2}+\sigma _{6}v_{3}+\sigma
_{7}v_{4}+\sigma _{5}\left( H^{\prime }+I^{\prime }\right) +\sigma
_{6}K^{\prime }+\sigma _{7}K-\sigma _{13}\}=0,$ }

\item {\small $\qquad j_{3}+\{\sigma _{6}v_{2}-\sigma _{8}v_{3}-\sigma
_{9}v_{4}-\sigma _{6}\left( H^{\prime }+I^{\prime }\right) -\sigma
_{8}K^{\prime }-\sigma _{9}K+\sigma _{15}\}=0,$ }

\item {\small $\qquad j_{4}+\{\sigma _{7}v_{2}-\sigma _{9}v_{3}-\sigma
_{10}v_{4}-\sigma _{7}\left( H^{\prime }+I^{\prime }\right) -\sigma
_{9}K^{\prime }-\sigma _{10}K+\sigma _{16}\}=0,$ }

\item {\small \qquad\ $k_{2}+\left\{ -\sigma _{5}w_{2}+\sigma
_{6}w_{3}+\sigma _{7}w_{4}-\sigma _{5}K^{\prime }+\sigma _{6}I^{\prime
}+\sigma _{7}J^{\prime }+\sigma _{19}\right\} =0,$ }

\item {\small $\qquad k_{3}+\left\{ \sigma _{6}w_{2}-\sigma _{8}w_{3}-\sigma
_{9}w_{4}+\sigma _{6}K^{\prime }-\sigma _{8}I^{\prime }-\sigma _{9}J^{\prime
}-\sigma _{21}\right\} =0,$ }

\item {\small $\qquad k_{4}+\left\{ \sigma _{7}w_{2}-\sigma _{9}w_{3}-\sigma
_{10}w_{4}+\sigma _{7}K^{\prime }-\sigma _{9}I^{\prime }-\sigma
_{10}J^{\prime }+\sigma _{19}\right\} =0.$ }
\end{description}

{\small Conversely, if (5.5) holds, then from (5.3), we get $\overline{H}%
_{\perp k}=\overline{I}_{\perp k}=\overline{J}_{\perp k}=\overline{K}_{\perp
k}=\overline{H}_{\perp k}^{\prime }=\overline{I^{\prime }}_{\perp k}=%
\overline{J^{\prime }}_{\perp k}=\overline{K^{\prime }}_{\perp k}=0$ and $%
\overline{h_{i}}=\overline{j_{i}}=\overline{k_{i}}=0.$ Thus the space $%
\left( M,\overline{L}\right) $ is a Berwald space. }

\begin{description}
\item[Case(vi)] {\small When $\sigma _{2}=0,\sigma _{3}$ $\neq 0$, $\sigma
_{4}=0.$ }
\end{description}

{\small In this case, if a four- dimensional Landsberg space is $\sigma -$%
conformally flat, then from the equations (5.2), (3.5) and (4.9), we see
that (5.4) reduces to }

\begin{description}
\item[$\left( 5.6\right) $] {\small $\frac{H;_{2}}{H;_{4}}=\frac{I;_{2}}{%
I;_{4}}=\frac{J;_{2}}{J;_{4}}=\frac{K;_{2}}{K;_{4}}=\frac{H^{\prime };_{2}}{%
H^{\prime };_{4}}=\frac{I^{\prime };_{2}}{I^{\prime };_{4}}=\frac{J^{\prime
};_{2}}{J^{\prime };_{4}}=\frac{K^{\prime };_{2}}{K^{\prime };_{4}}=\frac{%
\sigma _{7}}{\sigma _{5}}=\frac{\sigma _{9}}{\sigma _{6}}=\frac{\sigma _{10}%
}{\sigma _{7}}$ }
\end{description}

{\small and }

\begin{description}
\item {\small $\qquad h_{2}+\{-\sigma _{5}u_{2}+\sigma _{6}u_{3}+\sigma
_{7}u_{4}+\sigma _{5}\left( J+J^{\prime }\right) +\sigma _{6}I+\sigma
_{7}K^{\prime }-\sigma _{12}\}=0,$ }

\item {\small $\qquad h_{3}+\{\sigma _{6}u_{2}-\sigma _{8}u_{3}-\sigma
_{9}u_{4}-\sigma _{6}\left( J+J^{\prime }\right) -\sigma _{8}I-\sigma
_{9}K^{\prime }+\sigma _{14}\}=0,$ }

\item {\small $\qquad h_{4}+\{\sigma _{7}u_{2}-\sigma _{9}u_{3}-\sigma
_{10}u_{4}-\sigma _{7}\left( J+J^{\prime }\right) -\sigma _{9}I-\sigma
_{10}K^{\prime }+\sigma _{15}\}=0,$ }

\item {\small $\qquad j_{2}+\{-\sigma _{5}v_{2}+\sigma _{6}v_{3}+\sigma
_{7}v_{4}+\sigma _{5}\left( H^{\prime }+I^{\prime }\right) +\sigma
_{6}K^{\prime }+\sigma _{7}K-\sigma _{13}\}=0,$ }

\item {\small $\qquad j_{3}+\{\sigma _{6}v_{2}-\sigma _{8}v_{3}-\sigma
_{9}v_{4}-\sigma _{6}\left( H^{\prime }+I^{\prime }\right) -\sigma
_{8}K^{\prime }-\sigma _{9}K+\sigma _{15}\}=0,$ }

\item {\small $\qquad j_{4}+\{\sigma _{7}v_{2}-\sigma _{9}v_{3}-\sigma
_{10}v_{4}-\sigma _{7}\left( H^{\prime }+I^{\prime }\right) -\sigma
_{9}K^{\prime }-\sigma _{10}K+\sigma _{16}\}=0,$ }

\item {\small \ $\qquad k_{2}+\left\{ -\sigma _{5}w_{2}+\sigma
_{6}w_{3}+\sigma _{7}w_{4}-\sigma _{5}K^{\prime }+\sigma _{6}I^{\prime
}+\sigma _{7}J^{\prime }+\sigma _{19}\right\} =0,$ }

\item {\small $\qquad k_{3}+\left\{ \sigma _{6}w_{2}-\sigma _{8}w_{3}-\sigma
_{9}w_{4}+\sigma _{6}K^{\prime }-\sigma _{8}I^{\prime }-\sigma _{9}J^{\prime
}-\sigma _{21}\right\} =0,$ }

\item {\small $\qquad k_{4}+\left\{ \sigma _{7}w_{2}-\sigma _{9}w_{3}-\sigma
_{10}w_{4}+\sigma _{7}K^{\prime }-\sigma _{9}I^{\prime }-\sigma
_{10}J^{\prime }+\sigma _{19}\right\} =0.$ }
\end{description}

{\small Conversely, if (5.6) holds, then from (5.3), we get $\overline{H}%
_{\perp k}=\overline{I}_{\perp k}=\overline{J}_{\perp k}=\overline{K}_{\perp
k}=\overline{H}_{\perp k}^{\prime }=\overline{I^{\prime }}_{\perp k}=%
\overline{J^{\prime }}_{\perp k}=\overline{K^{\prime }}_{\perp k}=0$ and $%
\overline{h_{i}}=\overline{j_{i}}=\overline{k_{i}}=0.$ Thus the space $%
\left( M,\overline{L}\right) $ is a Berwald space. }

\begin{description}
\item[Case(vii)] {\small When $\sigma _{2}=0,\sigma _{3}$=$0$, $\sigma
_{4}\neq 0.$ }
\end{description}

{\small In this case, if a four- dimensional Landsberg space is $\sigma -$%
conformally flat, then from the equations (5.2), (3.5) and (4.10), we see
that (5.4) reduces to }

\begin{description}
\item[$\left( 5.7\right) $] {\small $\frac{H;_{2}}{H;_{3}}=\frac{I;_{2}}{%
I;_{3}}=\frac{J;_{2}}{J;_{3}}=\frac{K;_{2}}{K;_{3}}=\frac{H^{\prime };_{2}}{%
H^{\prime };_{3}}=\frac{I^{\prime };_{2}}{I^{\prime };_{3}}=\frac{J^{\prime
};_{2}}{J^{\prime };_{3}}=\frac{K^{\prime };_{2}}{K^{\prime };_{3}}=\frac{%
\sigma _{6}}{\sigma _{5}}=\frac{\sigma _{8}}{\sigma _{6}}=\frac{\sigma _{9}}{%
\sigma _{7}}$ }
\end{description}

{\small and }

\begin{description}
\item {\small $\qquad h_{2}+\{-\sigma _{5}u_{2}+\sigma _{6}u_{3}+\sigma
_{7}u_{4}+\sigma _{5}\left( J+J^{\prime }\right) +\sigma _{6}I+\sigma
_{7}K^{\prime }-\sigma _{12}\}=0,$ }

\item {\small $\qquad h_{3}+\{\sigma _{6}u_{2}-\sigma _{8}u_{3}-\sigma
_{9}u_{4}-\sigma _{6}\left( J+J^{\prime }\right) -\sigma _{8}I-\sigma
_{9}K^{\prime }+\sigma _{14}\}=0,$ }

\item {\small $\qquad h_{4}+\{\sigma _{7}u_{2}-\sigma _{9}u_{3}-\sigma
_{10}u_{4}-\sigma _{7}\left( J+J^{\prime }\right) -\sigma _{9}I-\sigma
_{10}K^{\prime }+\sigma _{15}\}=0,$ }

\item {\small $\qquad j_{2}+\{-\sigma _{5}v_{2}+\sigma _{6}v_{3}+\sigma
_{7}v_{4}+\sigma _{5}\left( H^{\prime }+I^{\prime }\right) +\sigma
_{6}K^{\prime }+\sigma _{7}K-\sigma _{13}\}=0,$ }

\item {\small $\qquad j_{3}+\{\sigma _{6}v_{2}-\sigma _{8}v_{3}-\sigma
_{9}v_{4}-\sigma _{6}\left( H^{\prime }+I^{\prime }\right) -\sigma
_{8}K^{\prime }-\sigma _{9}K+\sigma _{15}\}=0,$ }

\item {\small $\qquad j_{4}+\{\sigma _{7}v_{2}-\sigma _{9}v_{3}-\sigma
_{10}v_{4}-\sigma _{7}\left( H^{\prime }+I^{\prime }\right) -\sigma
_{9}K^{\prime }-\sigma _{10}K+\sigma _{16}\}=0,$ }

\item {\small \qquad\ $k_{2}+\left\{ -\sigma _{5}w_{2}+\sigma
_{6}w_{3}+\sigma _{7}w_{4}-\sigma _{5}K^{\prime }+\sigma _{6}I^{\prime
}+\sigma _{7}J^{\prime }+\sigma _{19}\right\} =0,$ }

\item {\small $\qquad k_{3}+\left\{ \sigma _{6}w_{2}-\sigma _{8}w_{3}-\sigma
_{9}w_{4}+\sigma _{6}K^{\prime }-\sigma _{8}I^{\prime }-\sigma _{9}J^{\prime
}-\sigma _{21}\right\} =0,$ }

\item {\small $\qquad k_{4}+\left\{ \sigma _{7}w_{2}-\sigma _{9}w_{3}-\sigma
_{10}w_{4}+\sigma _{7}K^{\prime }-\sigma _{9}I^{\prime }-\sigma
_{10}J^{\prime }+\sigma _{19}\right\} =0.$ }
\end{description}

{\small Conversely, if (5.7) holds, then from (5.3), we get $\overline{H}%
_{\perp k}=\overline{I}_{\perp k}=\overline{J}_{\perp k}=\overline{K}_{\perp
k}=\overline{H}_{\perp k}^{\prime }=\overline{I^{\prime }}_{\perp k}=%
\overline{J^{\prime }}_{\perp k}=\overline{K^{\prime }}_{\perp k}=0$ and $%
\overline{h_{i}}=\overline{j_{i}}=\overline{k_{i}}=0.$ Thus the space $%
\left( M,\overline{L}\right) $ is a Berwald space. Hence we have the
following: }

\begin{theorem}
{\small A four-dimensional conformally flat Landsberg space is a Berwald
space if and only if one of the following conditions is satisfied:}

\begin{description}
\item[1)] {\small \ the equations (5.4) hold for $\sigma _{2}\neq 0,\sigma
_{3}$ $\neq 0$, $\sigma _{4}\neq 0.$}

\item[2)] {\small $\sigma _{i}$ is orthogonal to $p^{i}$ and the equations
(5.4) hold for $\sigma _{2}\neq 0,\sigma _{3}$ $\neq 0$, $\sigma _{4}=0.$}

\item[3)] {\small $\sigma _{i}$ is orthogonal to $n^{i}$ and the equations
(5.4) hold for $\sigma _{2}\neq 0,\sigma _{3}$ $=0$, $\sigma _{4}\neq 0.$}

\item[4)] {\small \ $\sigma _{i}$ is orthogonal to $m^{i}$ and the equations
(5.4) hold for $\sigma _{2}=0,\sigma _{3}\neq 0$, $\sigma _{4}\neq 0.$}

\item[5)] {\small $\sigma _{i}$ is orthogonal to $n^{i}$ and $p^{i},$ and
the equations (5.5) hold for $\sigma _{2}\neq 0,\sigma _{3}$ $=0$, $\sigma
_{4}=0.$}

\item[6)] {\small $\sigma _{i}$ is orthogonal to $m^{i}$ and $p^{i},$ and
the equations (5.6) hold for $\sigma _{2}=0,\sigma _{3}$ $\neq 0$, $\sigma
_{4}=0.$}

\item[7)] {\small $\sigma _{i}$ is orthogonal to $m^{i}$ and $n^{i}, $ the
equations (5.7) hold for $\sigma _{2}=0,\sigma _{3}$=$0$, $\sigma _{4}\neq
0. $}
\end{description}
\end{theorem}
\noindent{\bf Acknowledgement}: The author is very much thankful to the referee for his valuable comments and suggestions towards the improvement of the paper.

\end{document}